\documentclass[12pt,a4paper,french,english]{article}
\usepackage[T1]{fontenc}
\usepackage[latin9]{inputenc}
\usepackage{color}
\usepackage{babel}
\makeatletter
\addto\extrasfrench{%
   \providecommand{\fg}{\ifdim\lastskip>\z@\unskip\fi~\frqq}%
}

\makeatother
\usepackage{amsthm}
\usepackage{amsmath}
\usepackage{amssymb}
\usepackage{graphicx}
\usepackage[unicode=true,pdfusetitle,
 bookmarks=true,bookmarksnumbered=false,bookmarksopen=false,
 breaklinks=false,pdfborder={0 0 0},backref=false,colorlinks=false]
 {hyperref}
\usepackage{breakurl}

\makeatletter

\numberwithin{equation}{section}
\numberwithin{figure}{section}
\theoremstyle{plain}
\newtheorem{thm}{\protect\theoremname}[section]
  \theoremstyle{definition}
  \newtheorem{defn}[thm]{\protect\definitionname}
  \theoremstyle{remark}
  \newtheorem{rem}[thm]{\protect\remarkname}
  \theoremstyle{plain}
  \newtheorem{prop}[thm]{\protect\propositionname}
  \theoremstyle{plain}
  \newtheorem{cor}[thm]{\protect\corollaryname}

\usepackage{a4wide}
\usepackage{color}
\usepackage{graphics}

\makeatother

  \addto\captionsenglish{\renewcommand{\corollaryname}{Corollary}}
  \addto\captionsenglish{\renewcommand{\definitionname}{Definition}}
  \addto\captionsenglish{\renewcommand{\propositionname}{Proposition}}
  \addto\captionsenglish{\renewcommand{\remarkname}{Remark}}
  \addto\captionsenglish{\renewcommand{\theoremname}{Theorem}}
  \addto\captionsfrench{\renewcommand{\corollaryname}{Corollaire}}
  \addto\captionsfrench{\renewcommand{\definitionname}{Définition}}
  \addto\captionsfrench{\renewcommand{\propositionname}{Proposition}}
  \addto\captionsfrench{\renewcommand{\remarkname}{Remarque}}
  \addto\captionsfrench{\renewcommand{\theoremname}{Théorème}}
  \providecommand{\corollaryname}{Corollary}
  \providecommand{\definitionname}{Definition}
  \providecommand{\propositionname}{Proposition}
  \providecommand{\remarkname}{Remark}
\providecommand{\theoremname}{Theorem}

\begin{document}

\title{Band structure of the Ruelle spectrum of contact Anosov flows}

\author{Frédéric Faure\\
 {\small{Institut Fourier, UMR 5582}}\\
 {\small{100 rue des Maths, BP74 }}\\
{\small{FR-38402 St Martin d'Hères, FRANCE}}\\
{\small{frederic.faure@ujf-grenoble.fr}}\\
{\small{http://www-fourier.ujf-grenoble.fr/\textasciitilde{}faure}}
\and\\
Masato Tsujii{\small{}}\\
{\small{Department of Mathematics, Kyushu University,}}\\
{\small{Moto-oka 744, Nishi-ku, Fukuoka, 819-0395, JAPAN }}\\
{\small{tsujii@math.kyushu-u.ac.jp}}\\
{\small{http://www2.math.kyushu-u.ac.jp/\textasciitilde{}tsujii }}}

\date{2013 May 5}
\maketitle
\begin{abstract}
If $X$ is a contact Anosov vector field on a smooth compact manifold
$M$ and $V\in C^{\infty}\left(M\right)$ it is known that the differential
operator $A=-X+V$ has some discrete spectrum called Ruelle-Pollicott
resonances in specific Sobolev spaces. We show that for $\left|\mathrm{Im}z\right|\rightarrow\infty$
the eigenvalues of $A$ are restricted to vertical bands and in the
gaps between the bands, the resolvent of $A$ is bounded uniformly
with respect to $\left|\mathrm{Im}\left(z\right)\right|$. In each
isolated band the density of eigenvalues is given by the Weyl law.
In the first band, most of the eigenvalues concentrate to the vertical
line $\mathrm{Re}\left(z\right)=\left\langle D\right\rangle _{M}$,
the space average of the function $D\left(x\right)=V\left(x\right)-\frac{1}{2}\mathrm{div}X_{\mid E_{u}\left(x\right)}$where
$E_{u}$ is the unstable distribution. This band spectrum gives an
asymptotic expansion for dynamical correlation functions.
\end{abstract}
\selectlanguage{french}%
~\newpage
\begin{abstract}
\textbf{Titre: ``Structure en bandes du spectre de Ruelle des flots
Anosov de contact''}.

Si $X$ est un champ de vecteur Anosov de contact sur une variété
compacte lisse $M$ et si $V\in C^{\infty}\left(M\right)$, il est
connu que l'opérateur différentiel $A=-X+V$ a du spectre discret
appelé résonances de Ruelle-Pollicott dans des espaces de Sobolev
spécifiques. On montre que pour $\left|\mathrm{Im}z\right|\rightarrow\infty$
les valeurs propres de $A$ sont inclues dans des bandes verticales
and que dans les gaps entre ces bandes la résolvante de $A$ est bornée
uniformément par rapport à $\left|\mathrm{Im}\left(z\right)\right|$.
Dans chaque bande isolée, la densité des valeurs propres est donnée
par une loi de Weyl. Dans la première bande, la plupart des valeurs
propres se concentrent sur la ligne verticale $\mathrm{Re}\left(z\right)=\left\langle D\right\rangle _{M}$,
qui est la moyenne spatiale de la fonction $D\left(x\right)=V\left(x\right)-\frac{1}{2}\mathrm{div}X_{\mid E_{u}\left(x\right)}$
où $E_{u}$ est la distribution instable. Ce spectre en bande permet
d'exprimer le comportement asymptotique des fonctions de correlations
dynamiques.
\end{abstract}
\selectlanguage{english}%

\section{Introduction}

In this paper we announce some results concerning the Ruelle-Pollicott
spectrum of transfer operators associated to contact Anosov flows
\cite{faure-tsujii_anosov_flows_13}. Let $X$ be a smooth vector
field on a compact manifold $M$ and suppose that $X$ generates a
contact Anosov flow.

The Ruelle-Pollicott spectrum of contact Anosov flows has been studied
since a long time due to its importance to describe the precise behavior
and decay of time correlation functions for large time. From this,
one can deduce fine statistical properties of the dynamics of the
flow such as exponential convergence towards equilibrium (i.e. mixing)
or central limit theorem for the Birkhoff average of functions. The
Ruelle-Pollicott spectrum is also useful to get some precise asymptotic
counting of periodic orbits.

Recent results show that the Ruelle-Pollicott resonances are the discrete
eigenvalues of the generator $\left(-X\right)$ seen as a differential
operator in some specific Sobolev spaces of distributions $\mathcal{H}\subset\mathcal{D}'\left(M\right)$
\cite{liverani_butterley_07,fred_flow_09,liverani_giulietti_2012}.
A more precise description of the structure of this spectrum has been
obtained in \cite{tsujii_08,tsujii_FBI_10} where it is shown that
in the asymptotic limit $\left|\mathrm{Im}z\right|\rightarrow\infty$
the spectrum is on the domain $\mathrm{Re}\left(z\right)\leq\gamma_{0}^{+}$
with some explicit ``gap'' $\gamma_{0}^{+}<0$ given below. More
generally these results can be extended to the operator $A=-X+V$
where $V\in C^{\infty}\left(M\right)$ is a smooth function called
``potential''.

In this paper we improve the description of the structure of this
Ruelle-Pollicott spectrum. The main results are stated in Theorem
\ref{thm:band-structure}. They show that the Ruelle-Pollicott spectrum
of the first order differential operator $A=-X+V$ has some band structure
in the asymptotic limit $\left|\mathrm{Im}z\right|\rightarrow\infty$,
i.e. is contained in the union of vertical bands $\mathbf{B}_{k}=\left\{ z\in\mathbb{C},\mathrm{Re}\left(z\right)\in\left[\gamma_{k}^{-},\gamma_{k}^{+}\right]\right\} $,
$k\geq0$ with $\gamma_{k+1}^{\pm}<\gamma_{k}^{\pm}$. The values
$\gamma_{k}^{+},\gamma_{k}^{-}$ are given explicitly in (\ref{eq:rk+-})
by the maximum (respect. minimum) of the time averaged along trajectories
of a function $D\in C^{\infty}\left(M\right)$ called ``damping function''
given by $D=V-\frac{1}{2}\mathrm{div}X_{\mid E_{u}}$. If the band
$\mathbf{B}_{k}$ is isolated from the others by an asymptotic spectral
gap (i.e. $\gamma_{k+1}^{+}<\gamma_{k}^{-}$) then the norm of resolvent
of $A$ is bounded in this gap uniformly with respect to $\left|\mathrm{Im}\left(z\right)\right|$.
Theorem \ref{thm:band-structure} shows that the spectrum in every
isolated band $\mathbf{B}_{k}$ satisfies a Weyl law, i.e. the number
$\mathcal{N}\left(b\right)$ of eigenvalues $z\in\mathbf{B}_{k}$
satisfying $\mathrm{Im}\left(z\right)\in\left[b,b+b^{\varepsilon}\right]$
is given by%
\footnote{The notation $\mathcal{N}\left(b\right)\asymp\left|b\right|^{d+\varepsilon}$
means that $\exists C>0$ independent of $b$ s.t. $\frac{1}{C}\left|b\right|^{d+\varepsilon}\leq\mathcal{N}\left(b\right)\leq C\left|b\right|^{d+\varepsilon}$.%
} $\mathcal{N}\left(b\right)/b^{\varepsilon}\asymp b^{d}$ as $b\rightarrow\infty$
for any $\varepsilon>0$, where $\mathrm{dim}M=2d+1$. The assumption
that the band is isolated is not needed for the upper bound. A better
result for the upper bound of this Weyl law is given in \cite{dyatlov_Ruelle_resonances_2012}:
it is shown  that for any radius $C_{0}>0$, the number of resonances
in the disk $D\left(ib,C_{0}\right)$ of center $ib$ is $O\left(b^{d}\right)$
(i.e. this is the case $\varepsilon=0$).

Concerning the most interesting ``external band'' $\mathbf{B}_{0}=\left\{ z\in\mathbb{C},\mathrm{Re}\left(z\right)\in\left[\gamma_{0}^{-},\gamma_{0}^{+}\right]\right\} $,
supposing that it is isolated ($\gamma_{1}^{+}<\gamma_{0}^{-}$),
it is shown in Theorem \ref{thm:concerntration} that most of the
resonances in the band $\mathbf{B}_{0}$ accumulate on the vertical
line $\mathrm{Re}\left(z\right)=\left\langle D\right\rangle _{M}$
given by the space average of the function $D$. This is due to ergodicity.
This problem is then closely related to the description of the spectrum
of the damped wave equation \cite{sjostrand_2000}. Finally Corollary
\ref{cor:correlation_decay} shows that dynamical correlation functions
can be expanded over the infinite spectrum contained in the first
band $\mathbf{B}_{0}$.

In the forthcoming paper \cite{faure-tsujii_anosov_flows_13} we will
consider the special case $V=V_{0}=\frac{1}{2}\mathrm{div}X_{\mid E_{u}}$
for which the damping function vanishes $D=0$, $\gamma_{0}^{\pm}=0$,
i.e. the Ruelle-Pollicott resonances of the external band accumulate
on the imaginary axis. However, $V_{0}$ is not smooth and this requires
an extension of the theory.

From Selberg theory and representation theory, this particular band
structure is known for a long time in the case of homogeneous hyperbolic
manifolds $M=\Gamma\backslash SO_{1,n}/SO_{n-1}\equiv\Gamma\backslash T_{1}^{*}\mathbb{H}^{n}$
where $\Gamma$ is a discrete co-compact subgroup. In that case, the
contact Anosov flow is the geodesic flow on the hyperbolic manifold
surface $\mathcal{N}=\Gamma\backslash\mathbb{H}^{n}=\Gamma\backslash SO_{1,n}/SO_{n}$.

Technically we use semiclassical analysis to study the spectrum of
the differential operator $A=-X+V$ \cite{martinez-01,zworski_book_2012}.
We consider the associated ``canonical dynamics'' in the phase space
$T^{*}M$ which is simply the lifted flow. The key observation is
that this canonical dynamics has a non-wandering set or ``trapped
set'' which is a smooth symplectic submanifold $K\subset T^{*}M$
and which is normally hyperbolic. This is the origin of the band structure
of the spectrum. The results presented in this paper have been already
obtained (among others) for a closely related problem, namely the
band structure of prequantum Anosov diffeomorphisms \cite{faure-tsujii_prequantum_maps_12}.
This approach has been originally developed on a simple model in \cite{fred-PreQ-06}.

In a recent paper \cite{dyatlov_resonance_band_2013}, Semyon Dyatlov
shows a band structure for resonances for a similar problem motivated
by scattering by black holes. The band structure he obtains also comes
from the property that the trapped set in his problem is symplectic
and normally hyperbolic but he assumes some smoothness for the (un)stable
foliations. One difficulty we have to deal with for Anosov flows is
the non smoothness of the (un)stable foliations. An other related
recent work is the paper of Nonnenmacher-Zworski \cite{nonenmacher_zworski_2013}
where they obtain Theorem \ref{thm:asymptotic-gap.-If} below but
for more general models including contact Anosov flow.

\selectlanguage{french}%
\newpage

\selectlanguage{english}%

\section{Contact Anosov flow}

%red 
\begin{center}{\color{red}\fbox{\color{black}\parbox{16cm}{
\begin{defn}
\label{def:Anosov}On a smooth Riemannian compact manifold $\left(M,g\right)$,
a smooth vector field $X$ generating a flow $\phi_{t}:M\rightarrow M,t\in\mathbb{R}$,
is \textbf{Anosov} (see Fig. \ref{fig:Anosov-flow.}), if there exists
an $\phi_{t}-$invariant decomposition of the tangent bundle $TM=E_{u}\oplus E_{s}\oplus E_{0},$
where $E_{0}=\mathbb{R}X$ and $C>0$,$\lambda>0$ such that for every
$t\geq0$ 
\begin{equation}
\left\Vert D\phi_{t/E_{s}}\right\Vert _{g}\leq Ce^{-\lambda t},\qquad\left\Vert D\phi_{-t/E_{u}}\right\Vert _{g}\leq Ce^{-\lambda t}.\label{eq:hyp_hyperbolic}
\end{equation}

\end{defn}
}}}\end{center}
\begin{rem}
In general the map $x\in M\rightarrow E_{u}\left(x\right),E_{s}\left(x\right)$
are only Hölder continuous. The ``structural stability theorem''
shows that Anosov vector field is a property robust under perturbation.
\end{rem}
\begin{center}
\begin{figure}
\begin{centering}
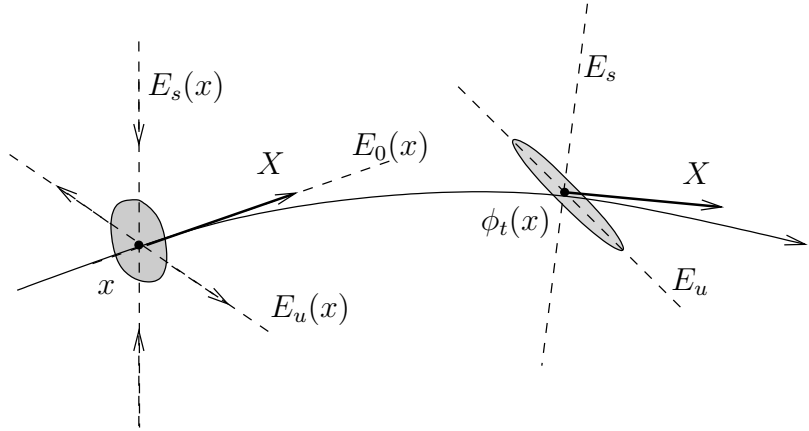
\par\end{centering}

\caption{\label{fig:Anosov-flow.}Anosov flow.}
\end{figure}

\par\end{center}

%red 
\begin{center}{\color{red}\fbox{\color{black}\parbox{16cm}{
\begin{defn}
The Anosov one form $\alpha\in C\left(T^{*}M\right)$ is defined by
$\mathrm{Ker}\alpha=E_{u}\oplus E_{s}$, $\alpha\left(X\right)=1$.
$X$ is a \textbf{contact Anosov vector field} if $\alpha$ is a smooth
contact one form i.e. $\left(d\alpha\right)_{\mid E_{u}\oplus E_{s}}$
is non degenerate (symplectic). 
\end{defn}
}}}\end{center}
\begin{rem}
If the case of a contact Anosov vector field we have that $\mathrm{dim}E_{u}=\mathrm{dim}E_{s}=d$,
with $\mathrm{dim}M=2d+1$ and that $dx=\alpha\wedge\left(d\alpha\right)^{d}$
is smooth volume form on $M$ preserved by the flow $\phi_{t}$.

As an example, the geodesic flow on a compact manifold $\mathcal{N}$
with negative sectional curvature (not necessary constant) defines
a contact Anosov flow on $T_{1}^{*}\mathcal{N}$. In that case the
Anosov one form $\alpha$ coincides with the canonical Liouville one
form $\xi dx$ on $T^{*}\mathcal{N}$.

We will assume that $X$ is a contact Anosov vector field on $M$
in the rest of this paper.
\end{rem}

\section{The transfer operator}

Let $V\in C^{\infty}\left(M\right)$ be a smooth function called ``potential''.

%red 
\begin{center}{\color{red}\fbox{\color{black}\parbox{16cm}{
\begin{defn}
The\textbf{ transfer operator} is the group of operators
\[
\hat{F}_{t}:\begin{cases}
C^{\infty}\left(M\right) & \rightarrow C^{\infty}\left(M\right)\\
v & \rightarrow e^{tA}v
\end{cases},\quad t\geq0
\]
with the generator
\[
A:=-X+V
\]
which is a first order differential operator.
\end{defn}
}}}\end{center}
\begin{rem}
~
\begin{itemize}
\item Since $X$ generates the flow $\phi_{t}$ we can write $\hat{F}_{t}v=\left(e^{\int_{0}^{t}V\circ\phi_{-s}ds}\right)v\left(\phi_{-t}\left(x\right)\right)$,
hence $\hat{F}_{t}$ acts as transport of functions by the flow with
multiplication by exponential of the function $V$ averaged along
the trajectory.
\item In the case $V=0$, the operator $\hat{F}_{t}$ is useful in order
to express ``dynamical correlation functions'' between $u,v\in C^{\infty}\left(M\right)$,
$t\in\mathbb{R}$: 
\begin{equation}
C_{u,v}\left(t\right):=\int_{M}u\cdot\left(v\circ\phi_{-t}\right)dx=\langle u,\hat{F}_{t}v\rangle_{L^{2}}\label{eq:time_correl}
\end{equation}
The study of these time correlation functions permits to establish
the mixing properties and other statistical properties of the dynamics
of the Anosov flow. In particular $u=cste$ is an obvious eigenfunction
of $A=-X$ with eigenvalue $z_{0}=0$. Since $\mathrm{div}X=0$ we
have that $\hat{F}_{t}$ is unitary in $L^{2}\left(M,dx\right)$ and
$iA=\left(iA\right)^{*}$is self-adjoint and has continuous spectrum
on the imaginary axis $\mathrm{Re}z=0$. In the next theorem we consider
more interesting functional spaces where the operator $A$ has discrete
spectrum but is non self-adjoint.
\end{itemize}
\end{rem}
%blue 
\begin{center}{\color{blue}\fbox{\color{black}\parbox{16cm}{
\begin{thm}[\cite{liverani_butterley_07}\cite{fred_flow_09}]
\label{thm:discrete_spectrum}''\textbf{discrete spectrum''}. If
$X$ is an Anosov vector field and $V\in C^{\infty}\left(M\right)$
then for every $C>0$, there exists a Hilbert space $\mathcal{H}_{C}$
with $C^{\infty}\left(M\right)\subset\mathcal{H}_{C}\subset\mathcal{D}'\left(M\right)$,
such that 
\[
A=-X+V\qquad:\mathcal{H}_{C}\rightarrow\mathcal{H}_{C}
\]
 has discrete spectrum on the domain $\mathrm{Re}\left(z\right)>-C\lambda$,
called \textbf{Ruelle-Pollicott resonances}, independent on the choice
of $\mathcal{H}_{C}$.
\end{thm}
}}}\end{center}
\begin{rem}
Concerning the meaning of these eigenvalues, notice that with the
choice $V=0$, if $\left(-X\right)v=zv$, $v$ is an invariant distribution
with eigenvalue $z=-a+ib\in\mathbb{C}$, then $v\circ\phi_{-t}=e^{-tX}v=e^{-at}e^{ibt}v$,
i.e. $a=-\mathrm{Re}\left(z\right)$ contributes as a damping factor
and $b=\mathrm{Im}\left(z\right)$ as a frequency in time correlation
function (\ref{eq:time_correl}). See corollary \ref{cor:correlation_decay}
below for a precise statement. Notice also the symmetry of the spectrum
under complex conjugation that$Av=zv$ implies $A\overline{v}=\overline{z}\overline{v}$.

We introduce now the following function that will play an important
role%
\footnote{Let $\mu_{g}$ be the induced Riemann volume form on $E_{u}\left(x\right)$
defined from the choice of a metric $g$ on $M$. As the usual definition
in differential geometry \cite[p.125]{taylor_tome1}, for tangent
vectors $u_{1},\ldots u_{d}\in E_{u}\left(x\right)$, $\mathrm{div}X_{\mid E_{u}}$
measures the rate of change of the volume of $E_{u}$ and is defined
by 
\[
\left(\mathrm{div}X_{\mid E_{u}}\left(x\right)\right)\cdot\mu_{g}\left(u_{1},\ldots u_{d}\right)=\lim_{t\rightarrow0}\frac{1}{t}\left(\mu_{g}\left(D\phi_{t}\left(u_{1}\right),\ldots,D\phi_{t}\left(u_{d}\right)\right)-\mu_{g}\left(u_{1},\ldots u_{d}\right)\right)
\]
Equivalently we can write that $\mathrm{div}X_{\mid E_{u}}\left(x\right)=\frac{d}{dt}\left(\mathrm{det}\left(D\phi_{t}\right)_{\mid E_{u}}\right)_{t=0}$.%
}
\begin{equation}
V_{0}\left(x\right):=\frac{1}{2}\mathrm{div}X_{\mid E_{u}}.\label{eq:def_V0}
\end{equation}
From (\ref{eq:hyp_hyperbolic}) we have $V_{0}\left(x\right)\geq\frac{1}{2}d\cdot\lambda$.
Since $E_{u}\left(x\right)$ is only Hölder in $x$ so is $V_{0}\left(x\right)$.
We will also consider the difference
\begin{equation}
D\left(x\right):=V\left(x\right)-V_{0}\left(x\right)\label{eq:def_D}
\end{equation}
and called it the ``\textbf{effective damping function}''. For simplicity
we will write:
\[
\left(\int_{0}^{t}D\right)\left(x\right):=\int_{0}^{t}\left(D\circ\phi_{-s}\right)\left(x\right)ds,\quad x\in M,
\]
for the Birkhoff average of $D$ along trajectories.
\end{rem}
%blue 
\begin{center}{\color{blue}\fbox{\color{black}\parbox{16cm}{
\begin{thm}[\cite{tsujii_08,tsujii_FBI_10}]
\label{thm:asymptotic-gap.-If}''\textbf{asymptotic gap''. }If
$X$ is a contact Anosov vector field on $M$ and $V\in C^{\infty}\left(M\right)$,
then for any $\varepsilon>0$ the Ruelle-Pollicott eigenvalues $\left(z_{j}\right)_{j}\in\mathbb{C}$
of $A=-X+V$ are contained in 
\[
\mathrm{Re}\left(z\right)\leq\gamma_{0}^{+}+\varepsilon
\]
up to finitely many exceptions and with 
\begin{equation}
\gamma_{0}^{+}=\lim_{t\rightarrow\infty}\sup_{x\in M}\frac{1}{t}\left(\int_{0}^{t}D\right)\left(x\right).\label{eq:gamma_0}
\end{equation}

\end{thm}
}}}\end{center}
\begin{rem}
See Figure \ref{fig:band_spectrum}(b). Notice that in the case $V=0$
we have $\gamma_{0}^{+}\leq-\frac{1}{2}d\cdot\lambda<0$.
\end{rem}

\section{Example of the geodesic flow on constant curvature surface}

A simple and well known example of contact Anosov flow is provided
by the geodesic flow on a surface $\mathcal{S}$ with constant negative
curvature. Precisely let $\Gamma<SL_{2}\mathbb{R}$ be a co-compact
discrete subgroup of $G=SL_{2}\mathbb{R}$ (i.e. such that $M:=\Gamma\backslash SL_{2}\mathbb{R}$
is compact). We suppose that $\left(-\mathrm{Id}\right)\in\Gamma$.
Then we have a natural identification that $M\equiv T_{1}^{*}\mathcal{S}$
is the unit cotangent bundle of the hyperbolic surface $\mathcal{S}:=\Gamma\backslash SL_{2}\mathbb{R}/SO_{2}\equiv\Gamma\backslash\mathbb{H}^{2}$.
Let $X$ be the left invariant vector field on $M$ given by the element
$X_{e}=\frac{1}{2}\left(\begin{array}{cc}
1 & 0\\
0 & -1
\end{array}\right)\in sl_{2}\mathbb{R}=T_{e}G$. Then $X$ is an Anosov contact vector field on $M$ and can be interpreted
as the geodesic flow on the surface $\mathcal{S}$. Using representation
theory, it is known that the Ruelle-Pollicott spectrum of the operator
$\left(-X\right)$ coincides with the zeros of the dynamical Fredholm
determinant. This dynamical Fredholm determinant is expressed as the
product of the Selberg zeta functions and gives the following result;
see figure \ref{fig:band_spectrum}(a). We refer to \cite{faure-tsujii_anosov_flows_13}
for further details.

%blue 
\begin{center}{\color{blue}\fbox{\color{black}\parbox{16cm}{
\begin{prop}
\label{thm:Ruelle spectrum_constant_curvature}If $X$ is the geodesic
flow on an hyperbolic surface $\mathcal{S}=\Gamma\backslash\mathbb{H}^{2}$
then the Ruelle-Pollicott eigenvalues $z$ of $\left(-X\right)$ are
of the form
\begin{equation}
z_{k,l}=-\frac{1}{2}-k\pm i\sqrt{\mu_{l}-\frac{1}{4}}\label{eq:z_mu_k}
\end{equation}
where $k\in\mathbb{N}$ and $\left(\mu_{l}\right)_{l\in\mathbb{N}}\in\mathbb{R}^{+}$
are the discrete eigenvalues of the hyperbolic Laplacian $\Delta$
on the surface $\mathcal{S}$. There are also $z_{n}=-n$ with $n\in\mathbb{N}^{*}$.
Each set $\left(z_{k,l}\right)_{l}$ with fixed $k$ will be called
the line $\mathbf{B}_{k}$. The ``Weyl law'' for $\Delta$ gives
the density of eigenvalues on each vertical line $\mathbf{B}_{k}$,
for $b\rightarrow\infty$,
\begin{equation}
\sharp\left\{ z_{k,l},\mbox{ }b<\mathrm{Im}\left(z_{k,l}\right)<b+1\right\} \asymp\left|b\right|\label{eq:Weyl_law_SL2R}
\end{equation}

\end{prop}
}}}\end{center}

\section{Band spectrum for general contact Anosov flow}

Proposition \ref{thm:Ruelle spectrum_constant_curvature} above shows
that the Ruelle-Pollicott spectrum for the geodesic flow on constant
negative surface has the structure of vertical lines $\mathbf{B}_{k}$
at $\mathrm{Re}z=-\frac{1}{2}-k$. In each line the eigenvalues are
in correspondence with the eigenvalues of the Laplacian $\Delta$.
We address now the question if this structure persists somehow for
geodesic flow on manifolds with negative (variable) sectional curvature
and more generally for any contact Anosov flow. In the next Theorem,
for a linear invertible map $L$, we note $\left\Vert L\right\Vert _{max}:=\left\Vert L\right\Vert $
and $\left\Vert L\right\Vert _{min}:=\left\Vert L^{-1}\right\Vert ^{-1}$. 

%blue 
\begin{center}{\color{blue}\fbox{\color{black}\parbox{16cm}{
\begin{thm}
\textbf{\label{thm:band-structure}}\cite{faure-tsujii_anosov_flows_13}\textbf{``asymptotic
band structure''.} If $X$ is a contact Anosov vector field on $M$
and $V\in C^{\infty}\left(M\right)$ then for every $C>0$, there
exists an Hilbert space $\mathcal{H}_{C}$ with $C^{\infty}\left(M\right)\subset\mathcal{H}_{C}\subset\mathcal{D}'\left(M\right)$,
such that for any $\varepsilon>0$, the Ruelle-Pollicott eigenvalues
$\left(z_{j}\right)_{j}\in\mathbb{C}$ of the operator $A=-X+V:\mathcal{H}_{C}\rightarrow\mathcal{H}_{C}$
on the domain $\mathrm{Re}\left(z\right)>-C\lambda$ are contained,
up to finitely many exceptions, in the union of finitely many bands
\[
z\in\bigcup_{k\geq0}\underbrace{\left[\gamma_{k}^{-}-\varepsilon,\gamma_{k}^{+}+\varepsilon\right]\times i\mathbb{R}}_{\mathrm{Band}\mbox{ }\mathbf{B}_{k}}
\]
with
\begin{align}
\gamma_{k}^{+} & =\lim_{t\rightarrow\infty}\left|\sup_{x}\frac{1}{t}\left(\left(\int_{0}^{t}D\right)\left(x\right)-k\log\left\Vert D\phi_{t}\left(x\right)_{/E_{u}}\right\Vert _{min}\right)\right|,\label{eq:rk+-}\\
\gamma_{k}^{-} & =\lim_{t\rightarrow\infty}\left|\inf_{x}\frac{1}{t}\left(\left(\int_{0}^{t}D\right)\left(x\right)-k\log\left\Vert D\phi_{t}\left(x\right)_{/E_{u}}\right\Vert _{max}\right)\right|
\end{align}
and where $D=V-V_{0}$ is the damping function (\ref{eq:def_D}).
In the gaps (i.e. between the bands) the norm of the resolvent is
controlled: there exists $c>0$ such that for every $z\notin\bigcup_{k\geq0}\mathbf{B}_{k}$
with $\left|\mathrm{Im}\left(z\right)\right|>c$ 
\[
\left\Vert \left(z-A\right)^{-1}\right\Vert \leq c.
\]
For some $k\geq0$, if the band $\mathbf{B}_{k}$ is ``isolated'',
i.e. $\gamma_{k+1}^{+}<\gamma_{k}^{-}$ and $\gamma_{k}^{+}<\gamma_{k-1}^{-}$
(this last condition is for $k\geq1$) then the number of resonances
in $\mathbf{B}_{k}$ obeys a ``\textbf{Weyl law}'': $\forall b>c$,
\begin{equation}
\frac{1}{c}\left|b\right|^{d}<\frac{1}{\left|b\right|^{\varepsilon}}\cdot\sharp\left\{ z_{j}\in\mathbf{B}_{k},b<\mathrm{Im}\left(z_{j}\right)<b+b^{\varepsilon}\right\} <c\left|b\right|^{d}\label{eq:Weyl_law}
\end{equation}
with $\dim M=2d+1$. The upper bound holds without the condition that
$\mathbf{B}_{k}$ is isolated.
\end{thm}
}}}\end{center}

\begin{center}
\begin{figure}
\begin{centering}
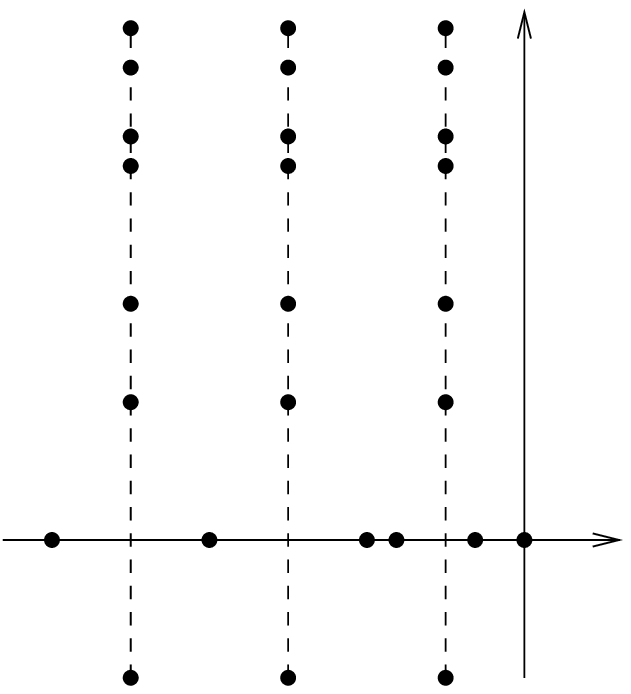$\qquad\quad$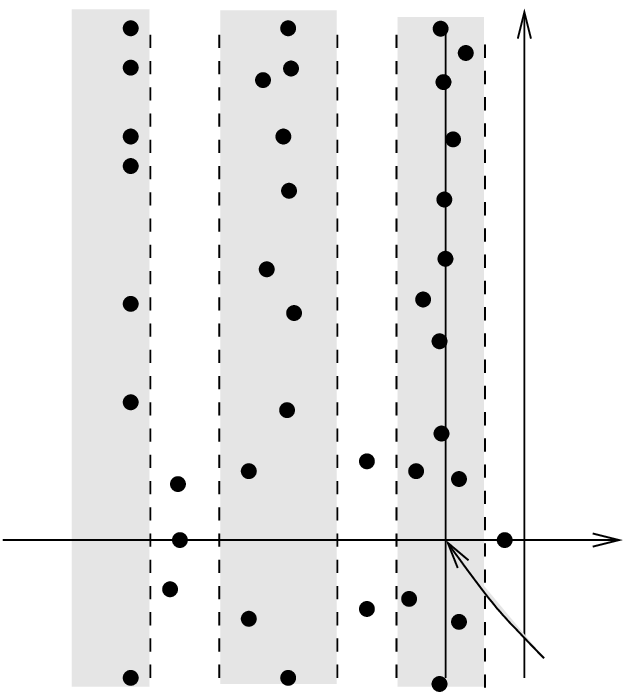
\par\end{centering}

\caption{\label{fig:band_spectrum}(a) For an hyperbolic surface $\mathcal{S}=\Gamma\backslash\mathbb{H}^{2}$,
the Ruelle-Pollicott spectrum of the geodesic vector field $-X$ given
by Proposition \ref{thm:Ruelle spectrum_constant_curvature}. It is
related to the eigenvalues of the Laplacian by (\ref{eq:z_mu_k}).
(b) For a general contact Anosov flow, the spectrum of $A=-X+V$ and
its asymptotic band structure given by Theorems \ref{thm:band-structure}
and \ref{thm:concerntration}.}
\end{figure}

\par\end{center}
\begin{rem}
We can compare Theorem \ref{thm:band-structure} with Proposition
\ref{thm:Ruelle spectrum_constant_curvature} in the special case
of the geodesic flow on a constant curvature surface $\mathcal{S}=\Gamma\backslash\mathbb{H}^{2}$:
we have $D\phi_{t}\left(x\right)_{/E_{u}}\equiv e^{t}$ hence $V_{0}=\frac{1}{2}$.
The choice of potential $V=0$ gives the constant damping function
$D=-\frac{1}{2}$, hence (\ref{eq:rk+-}) gives $\gamma_{k}^{+}=\gamma_{k}^{-}=-\frac{1}{2}-k$
as in Proposition \ref{thm:Ruelle spectrum_constant_curvature}.

In the forthcoming paper \cite{faure-tsujii_anosov_flows_13} we will
show that for a general contact Anosov vector field it is possible
to choose the potential $V=V_{0}$ (non smooth), giving $\gamma_{0}^{+}=\gamma_{0}^{-}=0$,
i.e. the first band is reduced to the imaginary axis and is isolated
from the second band by a gap, $\gamma_{1}^{+}<0$.
\end{rem}
%blue 
\begin{center}{\color{blue}\fbox{\color{black}\parbox{16cm}{
\begin{thm}
\label{thm:concerntration}\cite{faure-tsujii_anosov_flows_13}If
the external band $\mathbf{B}_{0}$ is isolated i.e. $\gamma_{1}^{+}<\gamma_{0}^{-}$,
then most of the resonances accumulate on the vertical line 
\[
\mathrm{Re}\left(z\right)=\left\langle D\right\rangle :=\frac{1}{\mathrm{Vol}\left(M\right)}\int_{M}D\left(x\right)dx
\]
in the precise sense that 
\begin{equation}
\frac{1}{\sharp\mathcal{B}_{b}}\sum_{z_{i}\in\mathcal{B}_{b}}\left|z_{i}-\left\langle D\right\rangle \right|\underset{b\rightarrow\infty}{\longrightarrow0},\qquad\mbox{with }\mathcal{B}_{b}:=\left\{ z_{i}\in\mathbf{B}_{0},\left|\mathrm{Im}\left(z_{i}\right)\right|<b\right\} .\label{eq:accumulation}
\end{equation}

\end{thm}
}}}\end{center}

\subsection*{Consequence for correlation functions expansion}

We mentioned the usefulness of dynamical correlation functions in
(\ref{eq:time_correl}). Let $\Pi_{j}$ denotes the finite rank spectral
projector associated to the eigenvalue $z_{j}$. The following Corollary
provides an expansion of correlation functions over the spectrum of
resonances of the first band $\mathbf{B}_{0}$. This is an infinite
sum.

%blue 
\begin{center}{\color{blue}\fbox{\color{black}\parbox{16cm}{
\begin{cor}
\label{cor:correlation_decay}Suppose that $\gamma_{1}^{+}<\gamma_{0}^{-}$.
Then for any $\varepsilon>0$, there exists $C_{\varepsilon}$, for
any $u,v\in C^{\infty}\left(M\right)$ and $t\geq0$,
\begin{equation}
\left|\langle u,\hat{F}_{t}v\rangle_{L^{2}}-\sum_{z_{j},\mathrm{Re}\left(z_{j}\right)\geq\gamma_{1}^{+}+\varepsilon}\langle u,\hat{F}_{t}\Pi_{j}v\rangle\right|\leq C_{\varepsilon}\left\Vert u\right\Vert _{\mathcal{H}_{C}'}\cdot\left\Vert u\right\Vert _{\mathcal{H}_{C}}\cdot e^{\left(\gamma_{1}^{+}+\varepsilon\right)t}.\label{eq:correlation_expansion}
\end{equation}
The infinite sum above converges fast because for arbitrary large
$m\geq0$ there exists $C_{m,\varepsilon}\left(u,v\right)\geq0$ such
that $\left|\langle u,\hat{F}_{t}\Pi_{j}v\rangle\right|\leq C_{m,\varepsilon}\left(u,v\right)\cdot e^{\left(\gamma_{0}^{+}+\varepsilon\right)t}\cdot\left|\mathrm{Im}\left(z_{j}\right)\right|^{-m}$
(except for a finite number of terms).
\end{cor}
}}}\end{center}

Eq.(\ref{eq:correlation_expansion}) is a refinement of decay of correlation
results of Dolgopyat \cite{dolgopyat_98}, Liverani \cite{liverani_contact_04},
Tsujii \cite[Cor.1.2]{tsujii_08,tsujii_FBI_10} and Nonnenmacher-Zworski
\cite[Cor.5]{nonenmacher_zworski_2013} where their expansion is a
finite sum over one or a finite number of leading resonances.

\subsection*{Outline of the proof}

The band structure and all related results presented above have already
been proven for the spectrum of Anosov prequantum map in \cite{faure-tsujii_prequantum_maps_12}.
An Anosov prequantum map $\tilde{f}:P\rightarrow P$ is an equivariant
lift of an Anosov diffeomorphism $f:M\rightarrow M$ on a principal
bundle $U\left(1\right)\rightarrow P\rightarrow M$ such that $\tilde{f}$
preserves a contact one form $\alpha$ (a connection on $P$). Therefore
$\tilde{f}:P\rightarrow P$ is very similar to the contact Anosov
flow $\phi_{t}:M\rightarrow M$ considered in this paper, that also
preserves a contact one form $\alpha$. Our proof of Theorem \ref{thm:band-structure}
is directly adapted from the proof given in \cite{faure-tsujii_prequantum_maps_12}.
We refer to this paper for more precisions on the proof and we use
the same notations below. The techniques rely on semiclassical analysis
adapted to the geometry of the contact Anosov flow lifted in the cotangent
space $T^{*}M$. In the limit $\left|\mathrm{Im}z\right|\rightarrow\infty$
of large frequencies under study, the semiclassical parameter is written
$\hbar:=1/\left|\mathrm{Im}z\right|$. We now sketch the main steps
of the proof.

\paragraph{Global geometrical description.}

$A=-X+V$ is a differential operator. Its principal symbol is the
function $\sigma\left(A\right)\left(x,\xi\right)=X_{x}\left(\xi\right)$
on phase space $T^{*}M$ (the cotangent bundle). It generates an Hamiltonian
flow which is simply the canonical lift of the flow $\phi_{t}$ on
$M$. Due to Anosov hypothesis on the flow in Definition \ref{def:Anosov},
the non-wandering set of the Hamiltonian flow is the continuous sub-bundle
$K=\mathbb{R}\alpha\subset T^{*}M$ where $\alpha$ is the Anosov
one form. $K$ is normally hyperbolic. This analysis has already been
used in \cite{fred_flow_09} for the semiclassical analysis of Anosov
flow (not necessary contact). With the additional hypothesis that
$\alpha$ is a smooth contact one form, this makes $K\backslash\left\{ 0\right\} $
a smooth symplectic submanifold of $T^{*}M$ (usually called the symplectization
of the contact one form $\alpha$) and normally hyperbolic. Let $\rho=\left(x,\xi\right)\in K$
be a point on the trapped set. Let $\hbar^{-1}=X_{x}\left(\xi\right)$
be its ``energy''. Let $\Omega=\sum_{j}dx^{j}\wedge d\xi^{j}$ be
the canonical symplectic form on $T^{*}M$ and consider the $\Omega$-orthogonal
splitting of the tangent space at $\rho\in K$: 
\begin{equation}
T_{\rho}\left(T^{*}M\right)=T_{\rho}K\oplus\left(T_{\rho}K\right)^{\perp}\label{eq:decomp}
\end{equation}
Due to hyperbolicity assumption, we have an additional decomposition
of the space 
\[
\left(T_{\rho}K\right)^{\perp}=E_{u}^{\left(2\right)}\oplus E_{s}^{\left(2\right)}
\]
 transverse to the trapped set into unstable/stable spaces.

\paragraph{Partition of unity.}

We decompose functions on the manifold using a microlocal partition
of unity of size $\hbar^{1/2-\varepsilon}$ with some $1/2>\varepsilon>0$,
that is refined as $\hbar\rightarrow0$. In each chart we use a canonical
change of variables adapted to the decomposition (\ref{eq:decomp}),
and construct an escape function adapted to the local splitting $E_{u}^{\left(2\right)}\oplus E_{s}^{\left(2\right)}$
above. This escape function has ``strong damping effect'' outside
a vicinity of size $O\left(\hbar^{1/2}\right)$ of the trapped set
$K$. We use this to define the anisotropic Sobolev space $\mathcal{H}_{C}$.
At the level of operators, we perform a decomposition similar to (\ref{eq:decomp})
and obtain a microlocal decomposition of the transfer operator $\hat{F}_{t}$
as a tensor product $\hat{F}_{t\mid T_{\rho}K}\otimes\hat{F}_{t\mid\left(T_{\rho}K\right)^{\perp}}$.
The first operator $\hat{F}_{t\mid T_{\rho}K}$ is unitary whereas
the second one $\hat{F}_{t\mid\left(T_{\rho}K\right)^{\perp}}$ has
discrete spectrum of resonances indexed by an integer $k\in\mathbb{N}$.
This is due to the choice of the escape function. We can construct
explicitly some approximate local spectral projectors $\Pi_{k}$ for
every value of $k$, and patching these locals expression together
we get global spectral operators for each band. The positions $\gamma_{k}^{\pm}$
of the band $\mathbf{B}_{k}$ come from estimates on the discrete
spectrum of the local operator $\hat{F}_{t\mid\left(T_{\rho}K\right)^{\perp}}$
restricted by the projector $\Pi_{k}$. We obtain results on the spectrum
of the generator $A$ from the results on the transfer operator $\hat{F}_{t}=e^{tA}$
by standard arguments.

The proof of the Weyl law is similar to the proof of J.Sjöstrand about
the damped wave equation \cite{sjostrand_2000} but needs more arguments.
The accumulation of resonances on the value $\left\langle D\right\rangle $
in Theorem \ref{thm:concerntration} given by the spatial average
of the damping function, Eq.(\ref{eq:accumulation}), uses the ergodicity
property of the Anosov flow and is also similar to the spectral results
obtained in \cite{sjostrand_2000} for the damped wave equation.

\selectlanguage{french}%
\bibliographystyle{plain}
\bibliography{/home/faure/articles/articles}
\selectlanguage{english}%

\end{document}